\documentclass[a4papet,10pt]{article}
\usepackage{amsfonts,amssymb,mathtools}
\usepackage{graphicx}
\usepackage{amsmath}
\usepackage[latin1]{inputenc} 
\usepackage{bbm}
\usepackage{pstricks}
\usepackage{psfrag}
\usepackage{appendix}
\usepackage{amsthm}
\usepackage{dsfont}
\usepackage{stackrel}
\usepackage{eufrak}
\usepackage{mathrsfs}
\usepackage{enumerate}
\usepackage{cite}
\usepackage{chngcntr}
\counterwithin{figure}{section}
\usepackage[a4paper,top=30mm,bottom=30mm,left=35mm,right=35mm]{geometry}
\numberwithin{equation}{section}

\newcommand{\probab}{\mathbb{P}}
\newcommand{\Ex}{\mathbb{E}}
\newcommand{\Rl}{\mathbb{R}}

\newtheorem{assum}{Assumption}[section]

\newtheoremstyle{mystyle}
  {}
  {}
  {\itshape}
  {}
  {\bfseries}
  {.}
  { }
  {\thmname{#1}\thmnumber{ #2}\thmnote{ (#3)}}

\theoremstyle{mystyle}
\newtheorem{theorem}{Theorem}[section]

\newtheorem{lemma}{Lemma}[section]
\newtheorem{proposition}{Proposition}[section]
\newtheorem{definition}{Definition}[section]

\title{Critical fluctuations in renewal models of statistical mechanics}

  \author{Marco Zamparo\footnote{Dipartimento di Fisica, Universit\`a degli Studi di Bari, via Amendola 173, Bari, I-70126, Italy
\newline \phantom{aaz} E-mail: \texttt{marco.zamparo@uniba.it}}}
  
\date{}

\begin{document}  
\maketitle

\begin{abstract}
  We investigate the sharp asymptotic behavior at criticality of the
  large fluctuations of extensive observables in renewal models of
  statistical mechanics, such as the Poland-Scheraga model of DNA
  denaturation, the Fisher-Felderhof model of fluids, the
  Wako-Sait\^o-Mu\~noz-Eaton model of protein folding, and the
  Tokar-Dreyss\'e model of strained epitaxy. These models amount to
  Gibbs changes of measure of a classical renewal process and can be
  identified with a constrained pinning model of polymers.  The
  extensive observables that enter the thermodynamic description turn
  out to be cumulative rewards corresponding to deterministic rewards
  that are uniquely determined by the waiting time and grow no faster
  than it. The probability decay with the system size of their
  fluctuations switches from exponential to subexponential at
  criticality, which is a regime corresponding to a discontinuous
  pinning-depinning phase transition. We describe such decay by
  proposing a precise large deviation principle under the assumption
  that the subexponential correction term to the waiting time
  distribution is regularly varying. This principle is in particular
  used to characterize the fluctuations of the number of renewals,
  which measures the DNA-bound monomers in the Poland-Scheraga model,
  the particles in the Fisher-Felderhof model and the Tokar-Dreyss\'e
  model, and the native peptide bonds in the
  Wako-Sait\^o-Mu\~noz-Eaton model.\\

\noindent Keywords: Renewal processes; Polymer pinning models;
Critical phenomena; Renewal-reward processes; Precise large
deviations; Regular varying tails\\

\noindent Mathematics Subject Classification 2020: Primary 60F10;
60K05; 60G50, Secondary 62E20; 60K35; 82B20; 82B23
\end{abstract}

\section{Introduction}
\label{sec:intro}

The Poland-Scheraga model of DNA denaturation \cite{PolSch1,PolSch2},
the Fisher-Felderhof model of fluids
\cite{FisFel1,FisFel2,FisFel3,FisFel4,Roepstorff}, the model of
protein folding introduced independently by Wako and Sait\^o first
\cite{WaSa1,WaSa2} and Mu\~noz and Eaton later
\cite{MuEa1,MuEa2,MuEa3}, and the Tokar-Dreyss\'e model of strained
epitaxy \cite{TokDre1,TokDre2,TokDre3} have been proved to share a
common regenerative structure in Ref.\ \cite{MarcoStat}, where they
have been named {\it renewal models of statistical mechanics}. In fact,
these models can be related to the pinning model of polymers
\cite{Giacomin,Frank}, which amounts to a Gibbs change of measure of a
classical discrete-time renewal process. Precisely, in
Ref.\ \cite{MarcoStat} they have been mathematically mapped into a
constrained pinning model obtained by the pinning model under the
condition that one of the renewals occurs at a predetermined time
corresponding to the system size. This constrained pinning model is
the prototype of renewal models of statistical mechanics and
introduces the language of renewal theory in the above lattice-gas
models of equilibrium statistical physics that apparently have nothing
to do with renewal theory.

The extensive variables that enter the thermodynamic description of
renewal models of statistical mechanics have been shown in
Ref.\ \cite{MarcoStat} to be cumulative rewards, supposing that each
renewal involves a deterministic reward that is uniquely determined by
the waiting time and grows no faster than it. Examples are the number
of DNA-bound monomers in the Poland-Scheraga model and the total
contact energy of a protein in the Wako-Sait\^o-Mu\~noz-Eaton model
\cite{MarcoStat}. The fluctuations of extensive variables in small
systems such as DNA molecules and proteins have become observable by
the recent advent of micromanipulation techniques \cite{Ritort},
motivating the development of a theoretical framework. A sharp large
deviation principle for multivariate extensive observables in renewal
models of statistical mechanics has been proposed in
Ref.\ \cite{MarcoStat}. This principle comes from more general sharp
large deviation principles established in Ref.\ \cite{MarcoProb} for
cumulative rewards in constrained and non-constrained pinning models
of polymers with broad-sense rewards taking values in a separable
Banach space.

Large deviation principles represent a general tool for describing the
exponential decay of probabilities of fluctuations in terms of a rate
function. A thorough study of the rate functions associated with
extensive observables in renewal models of statistical mechanics has
been provided in Ref.\ \cite{MarcoStat}.  There are however systems
for which the decay of probabilities is slower than exponential and
large deviation principles lose effectiveness due to a rate function
that displays a wide region of zeros. Such systems generalize a
polymer at the discontinuous pinning-depinning phase transition and
have been identified and named {\it critical} in
Ref.\ \cite{MarcoStat}.  Critical renewal models are found in the
theory of DNA and proteins since the DNA denaturation transition is
generally believed to be discontinuous \cite{Kafri1,Blossey1} and the
protein folding process is regarded as an ``all-or-none'' transition
\cite{Finkelstein1}, which is a microscopic analog of a discontinuous
phase transition.  Phase transitions are investigated by changing a
control parameter, such as the temperature or the denaturant
concentration, and the hallmark of discontinuous phase transitions is
a discontinuity in the graph of the expected value of some extensive
observable as a function of the control parameter. The number of
renewals, which counts for instance the pinned monomers in pinning
models of polymers and the DNA-bound monomers in the Poland-Scheraga
model, is the extensive observable commonly considered to characterize
the phase transition in renewal models of statistical mechanics
\cite{MarcoStat,Giacomin,Frank}. Criticality as defined in
Ref.\ \cite{MarcoStat} however does not require to introduce and
change a control parameter since the focus is on persistent
fluctuations that lead to subexponential decay of probabilities,
rather than on a comparison of the system in different
conditions. Despite the large amount of work devoted to pinning models
of polymers \cite{Giacomin,Frank}, which among others has elucidated
the sharp asymptotic behavior of partition functions by means of
methods from renewal theory, to the best of our knowledge criticality
in renewal models of statistical mechanics has never been investigated
in terms of fluctuations.

While a large deviation principle for extensive observables in renewal
models of statistical mechanics has been proposed in
Ref.\ \cite{MarcoStat} with no assumption on the waiting time
distribution, the study of critical fluctuations requires that some
hypothesis is made. This paper aims to describe the sharp asymptotics
of critical fluctuations, proposing a so-called precise large
deviation principle, under the assumption that the subexponential
correction term to the waiting time distribution is regularly
varying. This case covers pinning models of polymers and the
Poland-Scheraga model, for which the literature has mostly considered
polynomial-tailed corrections \cite{MarcoStat,Giacomin,Frank}.  The
paper is organized as follows.  Section \ref{sec:intro} introduces
pinning models and critical systems.  Section \ref{Results} presents
and discusses the main result of the paper and its application to the
number of renewals. The proof of this result is reported in section
\ref{proof}, which resorts to three appendices for the most technical
details.

\subsection{Pinning models}
\label{Models}

Let on a probability space $(\Omega,\mathcal{F},\probab)$ be given
independent and identically distributed random variables
$S_1,S_2,\ldots$ taking values in $\{1,2,\ldots\}\cup\{\infty\}$.  The
variable $S_i$ can be regarded as the \textit{waiting time} for the
$i$th occurrence at the \textit{renewal time} $T_i:=S_1+\cdots+S_i$ of
some event that is continuously renewed over time. The pinning model
considered in Refs.\ \cite{MarcoStat} and \cite{MarcoProb} makes use
of this formalism to describe a polymer that is pinned by a substrate
at the monomers $T_1,T_2,\ldots$. The polymer is supposed to consist
of $t\ge 1$ monomers so that the monomer $T_i$ contributes an energy
$-v(S_i)$ provided that $T_i\le t$, the real function $v$ defined on
$\{1,2,\ldots\}\cup\{\infty\}$ being called the
\textit{potential}. The state of the polymer is described by the law
$\probab_t$ on the measurable space $(\Omega,\mathcal{F})$ given by
the Gibbs change of measure
\begin{equation*}
\frac{d\probab_t}{d\probab}:=\frac{e^{H_t}}{Z_t},
\end{equation*}
where $H_t:=\sum_{i\ge 1}v(S_i)\mathds{1}_{\{T_i\le t\}}$ is the
\textit{Hamiltonian} and $Z_t:=\Ex[e^{H_t}]$ is the {\it partition
  function} ensuring normalization. The \textit{pinning model} is the
probabilistic model $(\Omega,\mathcal{F},\probab_t)$ supplied with the
hypotheses of aperiodicity and extensivity.  We say that the
\textit{waiting time distribution} $p:=\probab[S_1=\cdot\,]$ is
\textit{aperiodic} if its support $\mathcal{S}:=\{s\ge 1:p(s)>0\}$ is
nonempty and there does not exist an integer $\tau>1$ with the
property that $\mathcal{S}$ includes only some multiples of $\tau$. It
is worth observing that aperiodicity of $p$ can be obtained by simply
changing the time unit whenever $\probab[S_1<\infty]>0$.
\begin{assum}
\label{ass1}
The waiting time distribution $p$ is aperiodic.
\end{assum}
We say that the potential $v$ is \textit{extensive} if
$\limsup_{s\uparrow\infty}(1/s)\ln e^{v(s)}p(s)<+\infty$. Extensivity
is necessary to make the thermodynamic limit, as $t$ goes to infinity,
of the pinning model meaningful since
$Z_t\ge\Ex[e^{H_t}\mathds{1}_{\{S_1=t\}}]=e^{v(t)}p(t)$.
\begin{assum}
\label{ass2}
The potential $v$ is extensive.
\end{assum}

The \textit{constrained pinning model} where the last monomer is
forced to be always pinned by the substrate is the mathematical
skeleton of the Poland-Scheraga model, the Fisher-Felderhof model, the
Wako-Sait\^o-Mu\~noz-Eaton model, and the Tokar-Dreyss\'e model.
According to Refs.\ \cite{MarcoStat} and \cite{MarcoProb}, it
corresponds to the law $\probab_t^c$ on the measurable space
$(\Omega,\mathcal{F})$ defined through the change of measure
\begin{equation*}
\frac{d\probab_t^c}{d\probab}:=\frac{U_t e^{H_t}}{Z_t^c},
\end{equation*}
where $U_t:=\sum_{i\ge 1}\mathds{1}_{\{T_i=t\}}$ is the renewal
indicator taking value $1$ if and only if $t$ is a renewal and
$Z_t^c:=\Ex[U_te^{H_t}]$ is the partition function. As explained in
Ref.\ \cite{MarcoStat}, the constraint that $t$ is a renewal is needed
to set the system size, that is the number of monomers per strand for
the Poland-Scheraga model, the number of lattice sites for the
Fisher-Felderhof model and the Tokar-Dreyss\'e model, and the number
of peptide bonds for the Wako-Sait\^o-Mu\~noz-Eaton model, which turn
out to be $t$ (see \cite{MarcoStat}, section 3).  Aperiodicity of the
waiting time distribution gives $Z_t^c>0$ for all sufficiently large
$t$ (see \cite{MarcoProb}, section 1.1), thus ensuring that the
constrained pinning model is well-defined at least for such $t$.

\subsection{Deterministic rewards and critical systems}
\label{Known}

Let us suppose that the $i$th renewal involves a \textit{deterministic
  reward} $f(S_i)$, where $f$ is a function that maps
$\{1,2,\ldots\}\cup\{\infty\}$ in the Euclidean $d$-space $\Rl^d$.
The extensive observables of renewal models of statistical mechanics
are \textit{cumulative reward} by the time $t$ of the form
$W_t:=\sum_{i\ge 1}f(S_i)\mathds{1}_{\{T_i\le t\}}$ with $f(s)$ at
most of the order of magnitude of $s$ (see \cite{MarcoStat}, section
3). The number $N_t:=\sum_{i\ge 1}\mathds{1}_{\{T_i\le t\}}$ of
renewals by $t$ is the cumulative reward associated with the function
$f$ identically equal to 1. We aim to characterize the large
fluctuations of $W_t$ under the following hypothesis inherited from
Ref.\ \cite{MarcoStat}.
\begin{assum}
  \label{ass3}
If the support $\mathcal{S}$ of the waiting time distribution is
infinite, then $f(s)/s$ has a limit $r\in\Rl^d$ when $s$ goes to
infinity through $\mathcal{S}$.
\end{assum}

The study of the large fluctuations of $W_t$ on the exponential scale
has been carried out in Ref.\ \cite{MarcoStat}. This section collects
those results of Ref.\ \cite{MarcoStat} that introduce the problem of
critical fluctuations in renewal models of statistical mechanics and
that will serve the proof of the main contribution of the present
paper. In the thermodynamic limit, the scaled cumulative reward
$W_t/t$ converges in probability to a constant vector $\rho\in\Rl^d$
under assumptions \ref{ass1}, \ref{ass2}, and \ref{ass3}, which are
tacitly supposed to be satisfied in the sequel. In order to introduce
$\rho$, let us set $\ell:=\limsup_{s\uparrow\infty}(1/s)\ln
e^{v(s)}p(s)$, which fulfills $-\infty\le\ell<+\infty$ by assumption
\ref{ass2}, and if $\ell>-\infty$, then let us consider an ``effective
statistical weight'' $p_o$ for waiting times defined for all $s\ge 1$
by
\begin{equation}
  p_o(s):=e^{v(s)-\ell s}\,p(s).
\label{defpo}
\end{equation}
If $\ell=-\infty$ or $\ell>-\infty$ and $\sum_{s\ge 1}p_o(s)>1$, then
let $\zeta$ denote that unique real number larger than $\ell$ that
satisfies $\sum_{s\ge 1}e^{v(s)-\zeta s}\,p(s)=1$. Bearing in mind
that $\mathcal{S}$ is necessarily infinite when $\ell>-\infty$ and
letting $r$ be given by assumption \ref{ass3}, the vector $\rho$ turns
out to be
\begin{equation*}
  \rho:=\begin{cases}
  \frac{\sum_{s\ge 1} f(s)\,e^{v(s)-\zeta s}\,p(s)}{\sum_{s\ge 1} s\,e^{v(s)-\zeta s}\,p(s)} & \mbox{if } \ell=-\infty
  \mbox{ or }\ell>-\infty \mbox{ and } \sum_{s\ge 1}p_o(s)>1;\\
  \frac{\sum_{s\ge 1} f(s)\,p_o(s)}{\sum_{s\ge 1} s\,p_o(s)} & \mbox{if } \ell>-\infty, \mbox{ } \sum_{s\ge 1}p_o(s)=1,
  \mbox{ and } \sum_{s\ge 1}s\,p_o(s)<+\infty;\\
  r & \mbox{otherwise}.
  \end{cases}
\end{equation*}
We stress that $\sum_{s\ge 1} s\,e^{v(s)-\zeta s}\,p(s)$ is finite and
$\sum_{s\ge 1} f(s)\,e^{v(s)-\zeta s}\,p(s)$ exists due to assumption
\ref{ass3} whenever $\zeta$ is a real number larger than $\ell$.  The
following proposition states that $W_t/t$ convergences in probability
to $\rho$ as $t$ is sent to infinity (see \cite{MarcoStat}, theorem
4). It follows that $W_t/t$ convergences to $\rho$ also in mean since
$W_t/t$ takes bounded values with probability 1. Indeed, assumption
\ref{ass3} implies that there exists a positive constant $M<+\infty$
such that $\|f(s)\|\le Ms$ for all $s\in\mathcal{S}$, which gives
$\|W_t\|\le M\sum_{i\ge 1}S_i\mathds{1}_{\{T_i\le t\}}\le Mt$ with
probability 1.  Hereafter, $u\cdot v$ denotes the usual dot product
between $u$ and $v$ in $\Rl^d$ and $\|u\|:=\sqrt{u\cdot u}$ is the
Euclidean norm of $u$.
\begin{proposition}
\label{conv}
$\lim_{t\uparrow\infty}\probab_t^c[\|W_t/t-\rho\|\ge \delta]=0$ for any
$\delta>0$.
\end{proposition}
According to Ellis \cite{Ellisbook}, we say that $W_t/t$ {\it
  converges exponentially} to $\rho$ if for any $\delta>0$ there
exists a real number $\lambda>0$ such that
$\probab_t^c[\|W_t/t-\rho\|\ge\delta]\le e^{-\lambda t}$ for all
sufficiently large $t$. The following result improves proposition
\ref{conv} by identifying exponential convergence (see
\cite{MarcoStat}, theorem 4).
\begin{proposition}
\label{prop:conv}
$W_t/t$ converges exponentially to $\rho$ if and only if the conditions
  $\ell>-\infty$, $\sum_{s\ge 1}p_o(s)=1$, $\sum_{s\ge
    1}s\,p_o(s)<+\infty$, and $\rho\ne r$ are not simultaneously
  satisfied.
\end{proposition}
Proposition \ref{prop:conv} tells us that the convergence in
probability to $\rho$ of the scaled cumulative reward $W_t/t$ is
slower than exponential if $\ell>-\infty$, $\sum_{s\ge 1}p_o(s)=1$,
$\sum_{s\ge 1}s\,p_o(s)<+\infty$, and $\rho\ne r$. The facts that only
the condition $\rho\ne r$ involves the function $f$ and that such
condition is verified by most of $f$ justify the following definition
of critical model, which was originally proposed in
Ref.\ \cite{MarcoStat}. We stress that the condition $\sum_{s\ge
  1}p_o(s)=1$, if fulfilled, promotes $p_o$ to a probability
distribution on $\{1,2,\ldots\}$, which we call {\em effective waiting
  time distribution}.
\begin{definition}
The constrained pinning model is critical if $\ell>-\infty$,
$\sum_{s\ge 1}p_o(s)=1$, and $\sum_{s\ge 1}s\,p_o(s)<+\infty$.
\end{definition}
We point out that this notion of criticality is borrowed from the
literature on large deviation principles in statistical mechanics
\cite{El,Phys1}, where the focus is on the breaking of exponential
convergence. A different broader definition of critical system comes
from the theory of random polymers \cite{Giacomin,Frank}, where a
control parameter $\beta$ that plays the role of a binding energy can
drive a pinning-depinning phase transition. We recover the standard
framework of polymers by taking $v(s)=\beta$ for every $s$, and in
this framework our critical scenario corresponds to the conditions
$\ell=\limsup_{s\uparrow\infty}(1/s)\ln p(s)>-\infty$, $\sum_{s\ge
  1}s\,e^{-\ell s}\,p(s)<+\infty$, and $\beta=\beta_c$ with
$\beta_c:=-\ln\sum_{s\ge 1}e^{-\ell s}\,p(s)$. Clearly, $\beta_c$ is
finite if $\sum_{s\ge 1}s\,e^{-\ell s}\,p(s)<+\infty$.  These
conditions identify a discontinuous pinning-depinning phase transition
(see \cite{Giacomin}, theorem 2.1, or \cite{Frank}, theorem 7.4),
leaving out continuous phase transitions which preserve exponential
convergence. In order to clarify the point and make contact with the
literature \cite{MarcoStat,Giacomin,Frank} that looks at the expected
value of the number $N_t$ of renewals by $t$ to characterize the phase
transition, let us recall that $N_t$ is $W_t$ when the function $f$ is
identically equal to 1 and let us observe that this function satisfies
assumption \ref{ass3} with $r=0$. Then, $N_t/t$ converges in
probability and in mean to $\rho$ by proposition \ref{conv}. One can
verify (see \cite{MarcoStat}, section 4.4.1) that $\rho$ as a function
of $\beta$ is analytic on $(-\infty,+\infty)$ if $\ell=-\infty$ or
$\ell>-\infty$ and $\beta_c=-\infty$. If instead $\ell>-\infty$ and
$\beta_c>-\infty$, then $\rho$ is positive and analytic on the open
interval $(\beta_c,+\infty)$, continuous on the closed interval
$[\beta_c,+\infty)$, and equal to 0 for all $\beta<\beta_c$. In this
  case there is a phase transition for $\beta=\beta_c$. The phase
  transition is continuous, namely $\rho$ as a function of $\beta$ is
  continuous at $\beta_c$, if $\sum_{s\ge 1} s\,e^{-\ell
    s}\,p(s)=+\infty$. On the contrary, if $\sum_{s\ge 1} s\,e^{-\ell
    s}\,p(s)<+\infty$, then $\rho$ jumps from the value
  $w_c:=\frac{\sum_{s\ge 1} e^{-\ell s}\,p(s)}{\sum_{s\ge 1}
    s\,e^{-\ell s}\,p(s)}>0$ at $\beta_c$ to the value $0$ when
  $\beta$ is let to decrease and the phase transition is
  discontinuous.

Precise exponential rates for probability decays are provided by large
deviation principles. The cumulative reward $W_t$ satisfies a {\it
  large deviation principle} with {\it good rate function} according
to the following theorem (see \cite{MarcoStat}, theorem 1).
\begin{theorem}
\label{LDP}
There exists a proper convex lower semicontinuous function $I$ from
$\Rl^d$ to $[0,\infty]$ with compact level sets such that
\begin{enumerate}[{\upshape(a)}]
    \item $\liminf_{t\uparrow\infty}\frac{1}{t}\ln\probab_t^c\big[\frac{W_t}{t}\in G\big]\ge -\inf_{w\in G}\{I(w)\}$ 
    for each open set $G\subseteq\Rl^d$;
  \item $\limsup_{t\uparrow\infty}\frac{1}{t}\ln\probab_t^c\big[\frac{W_t}{t}\in F\big]\le-\inf_{w\in F}\{I(w)\}$
    for each Borel convex set $F\subseteq\Rl^d$ or closed set  $F\subseteq\Rl^d$.
\end{enumerate}
\end{theorem}
The function $I$ is the rate function and the compactness of its level
sets entails that $\inf_{w\in F}\{I(w)\}>0$ whenever $F$ is a closed
set that does not contain a zero of $I$ (see \cite{MarcoStat}, section
4.3).  The explicit expression of $I$ as the convex conjugate of the
scaled cumulant generating function has been given in
Ref.\ \cite{MarcoStat}, where the zeros of $I$ have in particular been
determined (see \cite{MarcoStat}, section 4.3) and the following
proposition has been deduced, $r$ being the vector of assumption
\ref{ass3}.
\begin{proposition}
Let $\mathcal{Z}$ be the set of zeroes of $I$. Then 
\begin{enumerate}[{\upshape(a)}]
  \item $\mathcal{Z}=\{\rho\}$ if the model is not critical;
  \item $\mathcal{Z}=\{(1-\alpha)r+\alpha \rho:\alpha\in[0,1]\}$ if the model is critical.
\end{enumerate}
\end{proposition}
The closed line segment $\mathcal{Z}$ is not a singleton at
criticality if $\rho\ne r$ and, according to the literature on large
deviation principles in statistical mechanics \cite{El,Phys1}, we call
it the {\it phase transition segment}.  Thus, whatever $f$ is, the
scaled cumulative reward $W_t/t$ converges exponentially to $\rho$ in
the non-critical scenario with exponential rate $\inf_{w\in
  F}\{I(w)\}>0$ for the probability of a fluctuation over a closed set
$F$ that does not contain $\rho$. On the contrary, convergence in
probability of $W_t/t$ to $\rho$ is slower than exponential in the
critical constrained pinning model provided that $f$ obeys $\rho\ne
r$. In this case, the above large deviation principle tells us that
the probability that $W_t/t$ fluctuates over a closed set $F$ that
does not intersect the phase transition segment $\mathcal{Z}$ decays
exponentially with rate $\inf_{w\in F}\{I(w)\}>0$, whereas it says
nothing about the fluctuations that reach $\mathcal{Z}$ and that
prevent exponential convergence. The situation $\rho=r$ constitutes an
exception because convergence is exponential even in the critical
scenario.

\section{Main results}
\label{Results}

The present paper describes, in a critical scenario with regularly
varying corrections to the waiting time distribution, the fluctuations
of the scaled cumulative reward $W_t/t$ that reach the phase
transition segment $\mathcal{Z}$.  Let us suppose that the constrained
pinning model is critical and that $\rho\ne r$.  We tackle the issue
by studying the fluctuations in the closed half-space
$\mathcal{H}_\alpha$ that contains the fraction $\alpha\in[0,1)$
  starting from $r$ of the segment $\mathcal{Z}$ and that is delimited
  by the hyperplane normal to it:
  \begin{equation*}
    \mathcal{H}_\alpha:=\Big\{w\in\Rl^d:[\rho-r]\cdot[w-(1-\alpha)r-\alpha\rho]\le 0\Big\}.
\end{equation*}
At large $t$, the probability
$\probab_t^c[W_t/t\in\mathcal{H}_\alpha]$ is dominated by the
fluctuations of $W_t/t$ over the phase transition segment that prevent
exponential convergence because we know that
$\probab_t^c[W_t/t\in\mathcal{H}_\alpha\cap F]$ decays exponentially
fast for every closed set $F$ that does not intersect $\mathcal{Z}$.

\subsection{A precise large deviation principle}

We determine the large-$t$ behavior of the probability
$\probab_t^c[W_t/t\in\mathcal{H}_\alpha]$ under the assumption that
there exist a real number $\kappa\ge 1$ and a slowly varying function
$\mathcal{L}$ such that for each $s\ge 1$
\begin{equation}
  p_o(s)=\frac{\mathcal{L}(s)}{s^{\kappa+1}},
\label{ass4}
\end{equation}
$p_o$ being the effective waiting time distribution defined by
(\ref{defpo}).  A measurable function $\mathcal{L}$ on the positive
semi-axis is {\it slowly varying} if it is positive on some
neighborhood of infinity and satisfies the scale-invariance property
$\lim_{x\uparrow+\infty}\mathcal{L}(\gamma x)/\mathcal{L}(x)=1$ for
any number $\gamma>0$. Trivially, a measurable function $\mathcal{L}$
with a positive limit at infinity is slowly varying. The simplest
non-trivial example is represented by the logarithm.  Due to slow
variation of $\mathcal{L}$, the sequence whose $s$th term is $p_o(s)$
is {\it regularly varying with index $-\kappa-1$}.  We refer to
\cite{RV} for the theory of slow and regular variation.  We point out
that the constraint $\kappa\ge 1$ reflects our focus on critical
models, which fulfill $\sum_{s\ge 1}s\,p_o(s)<+\infty$. Indeed, since
$\mathcal{L}(x)\ge x^{-\delta}$ for any fixed $\delta>0$ and all
sufficiently large $x$ (see \cite{RV}, proposition 1.3.6), $\kappa<1$
would entail $\sum_{s\ge 1}s\,p_o(s)=+\infty$.  The following theorem
is the main result of the paper.
\begin{theorem}
  \label{mainth}
Suppose that assumptions \ref{ass1}, \ref{ass2}, and \ref{ass3} are
satisfied, that the model is critical, and that $\rho\ne r$. If there
exist a real number $\kappa\ge 1$ and a slowly varying function
$\mathcal{L}$ such that (\ref{ass4}) holds, then for every
$\alpha\in[0,1)$
\begin{equation*}
  \lim_{t\uparrow\infty}\frac{\probab_t^c[W_t/t\in\mathcal{H}_\alpha]}{t^{1-\kappa}\mathcal{L}(t)}=\frac{1}{\sum_{s\ge 1}s\,p_o(s)}
  \begin{cases}
    \frac{\alpha}{1-\alpha}+\ln(1-\alpha) & \mbox{if } \kappa=1; \\
    \frac{1+(\alpha\kappa-1)(1-\alpha)^{-\kappa}}{\kappa(\kappa-1)} & \mbox{if } \kappa>1.
  \end{cases}
\end{equation*}
\end{theorem}

Theorem \ref{mainth} states a precise large deviation principle for
the cumulative reward $W_t$. It completes the picture of the large
fluctuations of extensive variables in renewal models of statistical
mechanics by supplying the leading order of probabilities at
criticality, where a large deviation principle loses effectiveness.
It is worth noting that the function $f$ of rewards enters only the
shape of the half-space $\mathcal{H}_\alpha$ through $\rho$ and $r$,
so that, basically, all extensive observables in critical models
display the same fluctuations, whose probability decays polynomially
in the system size under our assumption of regular variation.

The proof of theorem \ref{mainth} is provided in section \ref{proof}.
We stress that the cumulative reward is a random sum of random
variables, which can be written down explicitly by introducing the
number $N_t$ of renewals by $t$ as $W_t=\sum_{i=1}^{N_t}X_i$ with
$X_i:=f(S_i)$ for any $i$. Precise large deviations for random sums of
random variables with several types of subexponential distributions
have been investigated by many researchers
\cite{L1,L2,L3,L4,L5,L6,L7}. However, their work does not cover our
case in two respects. First, the hypothesis they have made is that the
random variables $X_1,X_2,\ldots$ are independent of the counting
process $N_t$, but this hypothesis is not satisfied by our
problem. Second, we had to implement the constraint that $t$ is a
renewal in order to set the system size of renewal models of
statistical mechanics at $t$, whereas they did not have such a special
need.  For these reasons, theorem \ref{mainth} is a new result that
requires a new proof, which is guided by the idea that the only
significant way in which a large fluctuation of $W_t$ occurs at
criticality is that one and only one of the waiting times takes a
large value. This idea underlies many heavy-tailed problems and is
known with the folklore name of {\em principle of a single big jump}
\cite{Foss}, whose a detailed picture can be given for sums of
independent and identically distributed random variables
\cite{Bigjump1,Bigjump2}.

\subsection{Critical fluctuations of $\boldsymbol{N_t}$}

As shown in Ref.\ \cite{MarcoStat}, renewal models of statistical
mechanics can be mapped in a standard constrained pinning model where
$v(s)=\beta$ for any $s$, $\beta$ being a real number that acts as a
control parameter. Precisely, $\beta$ is the binding energy for the
Poland-Scheraga model, the chemical potential for the Fisher-Felderhof
model and the Tokar-Dreyss\'e model, and an entropic loss for the
Wako-Sait\^o-Mu\~noz-Eaton model (see \cite{MarcoStat}, section
3). From section \ref{Known} we know that criticality is achieved if
$\ell=\limsup_{s\uparrow\infty}(1/s)\ln p(s)>-\infty$, $\sum_{s\ge
  1}s\,e^{-\ell s}\,p(s)<+\infty$, and $\beta=\beta_c$ with
$\beta_c:=-\ln\sum_{s\ge 1}e^{-\ell s}\,p(s)$. Let us suppose here
that $\ell>-\infty$ and that $\sum_{s\ge 1}s\,e^{-\ell
  s}\,p(s)<+\infty$ and let us investigate the fluctuations of the
number $N_t$ of renewals by $t$ for different values of $\beta$. The
extensive observable $N_t$ measures the DNA-bound monomers in the
Poland-Scheraga model, the particles in the Fisher-Felderhof model and
the Tokar-Dreyss\'e model, and the native peptide bonds in the
Wako-Sait\^o-Mu\~noz-Eaton model (see \cite{MarcoStat}, section 3).
By proposition \ref{conv}, $N_t/t$ converges in probability and in
mean to $\rho=\frac{\sum_{s\ge 1} e^{-\zeta s}\,p(s)}{\sum_{s\ge 1}
  s\,e^{-\zeta s}\,p(s)}$ for $\beta>\beta_c$ with $\zeta>\ell$
satisfying $\sum_{s\ge 1} e^{\beta-\zeta s}\,p(s)=1$, to
$\rho=w_c:=\frac{\sum_{s\ge 1} e^{-\ell s}\,p(s)}{\sum_{s\ge 1}
  s\,e^{-\ell s}\,p(s)}$ for $\beta=\beta_c$, and to $\rho=0$ for
$\beta<\beta_c$.  By combining parts (a) and (b) of theorem \ref{LDP}
we find for any $\beta$ and open interval $G$
\begin{equation*}
\lim_{t\uparrow\infty}\frac{1}{t}\ln\probab_t^c\bigg[\frac{N_t}{t}\in G\bigg]=-\inf_{w\in G}\big\{I(w)\big\}.
\end{equation*}
The rate function $I$ has been determined in Ref.\ \cite{MarcoStat} and
for each real $w$ takes the value (see \cite{MarcoStat}, section 4.4.2)
\begin{equation*}
  I(w):=\begin{cases}
  w(\beta_c-\beta)-\ell & \mbox{if } w\in[0,w_c];\\
  -w\ln\big[\sum_{s\ge 1}e^{\beta-\eta s}\,p(s)\big]-\eta & \mbox{if } w\in(w_c,1);\\
  -\ln\big[e^\beta p(1)\big] & \mbox{if } w=1;\\
  +\infty & \mbox{otherwise}
  \end{cases}
  +
  \begin{cases}
  \zeta & \mbox{if } \beta>\beta_c;\\
  \ell & \mbox{if } \beta\le\beta_c,
  \end{cases}
\end{equation*}
where $\eta$ is the unique number larger than $\ell$ that fulfills
$\frac{\sum_{s\ge 1} e^{-\eta s}\,p(s)}{\sum_{s\ge 1} s\,e^{-\eta
    s}\,p(s)}=w$.  The rate function $I$ has an affine stretch on
$[0,w_c]$, is analytic on $(w_c,1)$, and is continuously
differentiable at $w_c$ (see \cite{MarcoStat}, section 4.4.2). If
$\beta\ne\beta_c$, then the model is not critical and $I(w)=0$ only
for $w=\rho$. If instead $\beta=\beta_c$, then the model is critical
and $I(w)=0$ for all $w\in[0,w_c]$. The phase transition segment
$\mathcal{Z}$ is the closed interval $[0,w_c]$.

Theorem \ref{mainth} completes the study of the large fluctuations of
$N_t$ initiated in Ref.\ \cite{MarcoStat} by resolving for every
$\alpha\in[0,1)$ the subexponential decay at criticality of the
  probability $\probab_t^c[N_t\le\alpha w_c t]$. Indeed, under the
  assumption that there exist an index $\kappa\ge 1$ and a slowly
  varying function $\mathcal{L}$ such that $e^{-\ell
    s}\,p(s)=s^{-\kappa-1}\mathcal{L}(s)$ for all $s\ge 1$, as
  $\mathcal{H}_\alpha=(-\infty,\alpha w_c]$ theorem \ref{mainth} shows
that
\begin{equation*}
  \lim_{t\uparrow\infty}\frac{\probab_t^c[N_t\le\alpha w_c t]}{t^{1-\kappa}\mathcal{L}(t)}=\frac{1}{\sum_{s\ge 1}s\,e^{-\ell s}\,p(s)}
  \begin{cases}
    \frac{\alpha}{1-\alpha}+\ln(1-\alpha) & \mbox{if } \kappa=1; \\
    \frac{1+(\alpha\kappa-1)(1-\alpha)^{-\kappa}}{\kappa(\kappa-1)} & \mbox{if } \kappa>1.
  \end{cases}
\end{equation*}

\section{Proof of theorem \ref{mainth}}
\label{proof}

In this section we report the proof of theorem \ref{mainth},
postponing the most technical details in the appendices in order not
to interrupt the flow of the presentation. The proof is facilitated by
a natural change of measure. Regarding the effective waiting time
distribution $p_o$ as a new waiting time distribution, which is
non-defective because $\sum_{s\ge 1} p_o(s)=1$ by the hypothesis of
criticality, let us consider a new probability space
$(\Omega_o,\mathcal{F}_o,\probab_o)$ where a sequence $\{S_i\}_{i\ge
  1}$ of independent waiting times distributed according to $p_o$ is
given. Denoting by $\Ex_o$ the expectation with respect to
$\probab_o$, we have $\Ex_o[S_1]=\sum_{s\ge 1} s\,p_o(s)<+\infty$ and
$\Ex_o[f(S_1)]=\rho\,\Ex_o[S_1]$ again by criticality. The tail
probability $\mathcal{Q}$ of $S_1$ under the law $\probab_o$ is the
measurable function that maps each real number $x>0$ in
\begin{equation*}
\mathcal{Q}(x):=\probab_o[S_1>x]=\sum_{s\ge 1}\mathds{1}_{\{s>x\}} p_o(s).
\end{equation*}
The hypothesis that $p_o(s)=s^{-\kappa-1}\mathcal{L}(s)$ for all $s\ge
1$ with $\mathcal{L}$ a slowly varying function results in regular
variation with index $-\kappa$ for $\mathcal{Q}$, meaning that
$\lim_{x\uparrow+\infty}\mathcal{Q}(\gamma
x)/\mathcal{Q}(x)=\gamma^{-\kappa}$ for any number $\gamma>0$. Indeed,
the following lemma holds as demonstrated in appendix \ref{appendix1}.
\begin{lemma}
\label{lemmaQtail}
$\mathcal{Q}$ varies regularly with index $-\kappa$ and
$\lim_{x\uparrow+\infty}x^\kappa\mathcal{Q}(x)/\mathcal{L}(x)=1/\kappa$.
\end{lemma}

The transition to the new probabilistic model
$(\Omega_o,\mathcal{F}_o,\probab_o)$ proceeds as follows.  By
observing that $\prod_{i=1}^n e^{v(s_i)}p(s_i)=e^{\ell t}\prod_{i=1}^n
p_o(s_i)$ whenever $s_1+\cdots+s_n=t$, for every integer $t\ge 1$ and
Borel set $\mathcal{B}\subseteq\Rl^d$ we find
\begin{align}
\nonumber
Z_t^c\cdot\,\probab_t^c\bigg[\frac{W_t}{t}\in \mathcal{B}\bigg]&=\Ex\bigg[\mathds{1}_{\big\{\frac{W_t}{t}\in \mathcal{B}\big\}} U_te^{H_t}\bigg]\\
\nonumber
&=\sum_{n\ge 1}\Ex\bigg[\mathds{1}_{\big\{\frac{1}{t}\sum_{i=1}^n f(S_i)\in \mathcal{B}\big\}} \mathds{1}_{\{T_n=t\}}e^{\sum_{i=1}^n v(S_i)}\bigg]\\
\nonumber
&=\sum_{n\ge 1}\sum_{s_1\ge 1}\cdots\sum_{s_n\ge 1}\mathds{1}_{\big\{\frac{1}{t}\sum_{i=1}^n f(s_i)\in \mathcal{B}\big\}} \mathds{1}_{\{s_1+\cdots+s_n=t\}}
\prod_{i=1}^n e^{v(s_i)} p(s_i)\\
\nonumber
&=e^{\ell t}\sum_{n\ge 1}\sum_{s_1\ge 1}\cdots\sum_{s_n\ge 1}\mathds{1}_{\big\{\frac{1}{t}\sum_{i=1}^n f(s_i)\in \mathcal{B}\big\}} \mathds{1}_{\{s_1+\cdots+s_n=t\}}
\prod_{i=1}^n p_o(s_i)\\
&=e^{\ell t}\,\Ex_o\bigg[\mathds{1}_{\big\{\frac{W_t}{t}\in \mathcal{B}\big\}} U_t\bigg].
\label{primaident}
\end{align}
This identity with $\mathcal{B}=\Rl^d$ yields $Z_t^c=e^{\ell
  t}\,\Ex_o[U_t]$, which allows us to recast (\ref{primaident}) as
\begin{equation}
\label{mod2}
\Ex_o[U_t]\cdot\probab_t^c\bigg[\frac{W_t}{t}\in\mathcal{B}\bigg]=\Ex_o\bigg[\mathds{1}_{\big\{\frac{W_t}{t}\in \mathcal{B}\big\}} U_t\bigg].
\end{equation}
Formula (\ref{mod2}) is precisely the bridge between constrained
pinning models with respect to $(\Omega,\mathcal{F},\probab)$ and
$(\Omega_o,\mathcal{F}_o,\probab_o)$. We observe that
$\lim_{t\uparrow\infty}\Ex_o[U_t]=1/\Ex_o[S_1]$ since $p_o$ is
non-defective. This limit is due to an application of the renewal
theorem (see \cite{Feller1}, theorem 1 in chapter XIII.10) to the
renewal equation $\Ex_o[U_t]=\sum_{s=1}^tp_o(s)\,\Ex_o[U_{t-s}]$ valid
for every $t\ge 1$. The renewal equation is deduced by conditioning on
$T_1=S_1$ and then by using the fact that a renewal process starts
over at every renewal.

Fix $\alpha\in[0,1)$. If $t$ is a renewal, then $\sum_{i\ge
    1}S_i\mathds{1}_{\{T_i\le t\}}=t$ and the condition
  $W_t/t\in\mathcal{H}_\alpha$, namely
  $[\rho-r]\cdot[W_t/t-(1-\alpha)r-\alpha\rho]\le 0$, becomes
  $\sum_{i\ge 1}g(S_i)\mathds{1}_{\{T_i\le t\}}\le\alpha t$ with
  \begin{equation*}
    g(s):=\frac{[\rho-r]\cdot[f(s)-rs]}{\|\rho-r\|^2}.
\end{equation*}
This way, by setting $\mathcal{B}=\mathcal{H}_\alpha$ in (\ref{mod2})
we find for each $t\ge 1$
\begin{equation*}
  \Ex_o[U_t]\cdot\probab_t^c\bigg[\frac{W_t}{t}\in\mathcal{H}_\alpha\bigg]=
  \Ex_o\bigg[\mathds{1}_{\big\{\sum_{i\ge 1}g(S_i)\mathds{1}_{\{T_i\le t\}}\le\alpha t\big\}}U_t\bigg]=:\mathcal{E}(t).
\end{equation*}
Since $\lim_{t\uparrow\infty}\Ex_o[U_t]=1/\Ex_o[S_1]$ and
$\lim_{t\uparrow\infty}t^\kappa\mathcal{Q}(t)/\mathcal{L}(t)=1/\kappa$
by lemma \ref{lemmaQtail}, it follows that in order to prove theorem
\ref{mainth} it suffices to demonstrate that
\begin{align}
  \lim_{t\uparrow\infty}\frac{\mathcal{E}(t)}{t\,\mathcal{Q}(t)}
 \nonumber
  &=\frac{1}{\Ex_o[S_1]^2}\begin{cases}
    \frac{\alpha}{1-\alpha}+\ln(1-\alpha) & \mbox{if } \kappa=1;\\
    \frac{1+(\alpha\kappa-1)(1-\alpha)^{-\kappa}}{\kappa-1} & \mbox{if } \kappa>1
 \end{cases}\\
 &=\frac{1}{\Ex_o[S_1]^2}\bigg[\alpha(1-\alpha)^{-\kappa}-\int_{1-\alpha}^1\frac{dx}{x^\kappa}\bigg].
 \label{mod3}
\end{align}
We shall demonstrate (\ref{mod3}) by verifying a lower bound first and
an upper bound later. The features of the function $g$ we will use to
this aim are the equality $\Ex_o[g(S_1)]=\Ex_o[S_1]$ and the fact that
$g(s)/s$ goes to zero when $s$ is sent to infinity through
$\mathcal{S}$, which is the support of both $p$ and $p_o$. Since the
values $g(s)$ when $s\notin\mathcal{S}$ do not affect the problem, in
order to simplify the notations we redefine the function $g$ by
setting $g(s):=0$ for all $s\notin\mathcal{S}$ so that the limit
$\lim_{s\uparrow\infty}g(s)/s=0$ is valid.

\subsection{A lower bound}

In order to prove (\ref{mod3}) we show at first that
\begin{equation}
   \label{binf}
  \liminf_{t\uparrow\infty}\frac{\mathcal{E}(t)}{t\,\mathcal{Q}(t)}\ge
  \frac{1}{\Ex_o[S_1]^2}\bigg[\alpha(1-\alpha)^{-\kappa}-\int_{1-\alpha}^1\frac{dx}{x^\kappa}\bigg].
\end{equation}
This bound is trivial if $\alpha=0$ since the l.h.s.\ of (\ref{binf})
is not negative. In the case $\alpha>0$, it follows if we demonstrate
that for all real numbers $\gamma$, $\eta$, and $\epsilon$ such that
$1-\alpha<\gamma<\eta<1$ and $\epsilon>0$
\begin{equation}
\label{binf1}
  \liminf_{t\uparrow\infty}\frac{\mathcal{E}(t)}{t\,\mathcal{Q}(t)}\ge
  \frac{1-2\epsilon}{\Ex_o[S_1]^2}\bigg[(1-\gamma)\gamma^{-\kappa}-(1-\eta)\eta^{-\kappa}-\int_\gamma^\eta \frac{dx}{x^\kappa}\bigg].
\end{equation}
Indeed, (\ref{binf1}) gives (\ref{binf}) when $\gamma$ is sent to
$1-\alpha$, $\eta$ is sent to 1, and $\epsilon$ is sent to 0.

Let us assume $\alpha>0$ and let us pick $\gamma$, $\eta$, and
$\epsilon$ such that $1-\alpha<\gamma<\eta<1$ and $\epsilon>0$.  To
get at a proof of (\ref{binf1}) we observe that $3y/2-y^2/2\le 1$ for
any integer $y\in\mathbb{Z}$, so that for each $n\ge 2$ we find
\begin{equation*}
1\ge\frac{3}{2}\sum_{j=1}^n\mathds{1}_{\{S_j>\gamma t\}}-\frac{1}{2}\Bigg[\sum_{j=1}^n\mathds{1}_{\{S_j>\gamma t\}}\Bigg]^2
=\sum_{j=1}^n\mathds{1}_{\{S_j>\gamma t\}}-\sum_{j=1}^{n-1}\sum_{k=j+1}^n\mathds{1}_{\{S_j>\gamma t\}}\mathds{1}_{\{S_k>\gamma t\}}.
\end{equation*}
It follows that for all $t\ge 2$
\begin{align}
  \nonumber
  \mathcal{E}(t)&:=\Ex_o\bigg[\mathds{1}_{\big\{\sum_{i\ge 1}g(S_i)\mathds{1}_{\{T_i\le t\}}\le\alpha t\big\}}U_t\bigg]
  =\sum_{n=1}^t\Ex_o\bigg[\mathds{1}_{\big\{\sum_{i=1}^ng(S_i)\le\alpha t\big\}}\mathds{1}_{\{T_n=t\}}\bigg]\\
  \nonumber
  &\ge\sum_{n=1}^t\sum_{j=1}^n\Ex_o\bigg[\mathds{1}_{\big\{\sum_{i=1}^ng(S_i)\le\alpha t\big\}}\mathds{1}_{\{S_j>\gamma t\}}
    \mathds{1}_{\{S_1+\cdots + S_n=t\}}\bigg]+\\
  \nonumber
  &-\sum_{n=2}^t\sum_{j=1}^{n-1}\sum_{k=j+1}^n
  \Ex_o\bigg[\mathds{1}_{\big\{\sum_{i=1}^ng(S_i)\le\alpha t\big\}}\mathds{1}_{\{S_j>\gamma t\}}\mathds{1}_{\{S_k>\gamma t\}}
    \mathds{1}_{\{S_1+\cdots+S_n=t\}}\bigg]\\
  \nonumber
  &\ge\sum_{n=1}^tn\,\Ex_o\bigg[\mathds{1}_{\big\{\sum_{i=1}^ng(S_i)\le\alpha t\big\}}\mathds{1}_{\{S_n>\gamma t\}}\mathds{1}_{\{T_n=t\}}\bigg]
  -t^2\,\Ex_o\Big[\mathds{1}_{\{S_1>\gamma t\}}\mathds{1}_{\{S_2>\gamma t\}}U_t\Big]\\
  &\ge\sum_{n=1}^tn\,\Ex_o\bigg[\mathds{1}_{\big\{\sum_{i=1}^ng(S_i)\le\alpha t\big\}}\mathds{1}_{\{S_n>\gamma t\}}\mathds{1}_{\{T_n=t\}}\bigg]
  -\big[t\,\mathcal{Q}(\gamma t)\big]^2.
  \label{aux1}
\end{align}
Since $\alpha+\gamma-1>0$ by hypothesis, the limit
$\lim_{s\uparrow\infty}g(s)/s=0$ entails that there exists a positive
integer $t_1$ such that $g(s)\le (\alpha+\gamma-1)s$ whenever
$s>\gamma t_1$. If $t>t_1$, then the conditions $s>\gamma t$ and
$\sum_{i=1}^ng(S_i)<(\alpha+\gamma)(t-s)$ together imply
\begin{equation*}
\sum_{i=1}^ng(S_i)+g(s)\le(\alpha+\gamma)(t-s)+(\alpha+\gamma-1)s=(\alpha+\gamma)t-s\le\alpha t.
\end{equation*}
This way, from (\ref{aux1}) we get for every $t>t_1\ge 1$
\begin{align}
  \nonumber
  \mathcal{E}(t)&\ge\sum_{n=2}^t\sum_{s=1}^{t-n+1} n\,\Ex_o\bigg[\mathds{1}_{\big\{\sum_{i=1}^{n-1}g(S_i)+g(s)\le\alpha t\big\}}
    \mathds{1}_{\{T_{n-1}=t-s\}}\bigg]\,\mathds{1}_{\{s>\gamma t\}}\,p_o(s)-\big[t\,\mathcal{Q}(\gamma t)\big]^2\\
  \nonumber
  &\ge\sum_{s=1}^{t-1}\sum_{n=1}^{t-s} n\,\Ex_o\bigg[\mathds{1}_{\big\{\sum_{i=1}^ng(S_i)+g(s)\le\alpha t\big\}}
    \mathds{1}_{\{T_n=t-s\}}\bigg]\,\mathds{1}_{\{s>\gamma t\}}\,p_o(s)-\big[t\,\mathcal{Q}(\gamma t)\big]^2\\
  \nonumber
  &\ge\sum_{s=1}^{t-1}\sum_{n=1}^{t-s} n\,\Ex_o\bigg[\mathds{1}_{\big\{\sum_{i=1}^ng(S_i)< (\alpha+\gamma)(t-s)\big\}}
    \mathds{1}_{\{T_n=t-s\}}\bigg]\,\mathds{1}_{\{s>\gamma t\}}\,p_o(s)-\big[t\,\mathcal{Q}(\gamma t)\big]^2\\
  \nonumber
  &=\sum_{s=1}^{t-1}\Ex_o\bigg[N_{t-s}\mathds{1}_{\big\{\sum_{i\ge 1}g(S_i)\mathds{1}_{\{T_i\le t-s\}}< (\alpha+\gamma)(t-s)\big\}}U_{t-s}\bigg]\,
  \mathds{1}_{\{s>\gamma t\}}\,p_o(s)-\big[t\,\mathcal{Q}(\gamma t)\big]^2.
\end{align}
By introducing the restriction $s<\eta t$, for any $t>t_1$ we obtain
the further lower bound
\begin{align}
  \nonumber
 \mathcal{E}(t)&\ge\sum_{s\ge 1}\Ex_o\bigg[N_{t-s}\mathds{1}_{\big\{\sum_{i\ge 1}g(S_i)\mathds{1}_{\{T_i\le t-s\}}<(\alpha+\gamma)(t-s)\big\}}U_{t-s}\bigg]\,
  \mathds{1}_{\{\gamma t<s\le\eta t\}}\,p_o(s)-\big[t\,\mathcal{Q}(\gamma t)\big]^2\\
  \nonumber
  &\ge\sum_{s\ge 1}\Ex_o\big[N_{t-s}U_{t-s}\big]\,\mathds{1}_{\{\gamma t<s\le\eta t\}}\,p_o(s)-\big[t\,\mathcal{Q}(\gamma t)\big]^2+\\
  \label{aux2}
  &-\sum_{s\ge 1}(t-s)\,\Ex_o\bigg[\mathds{1}_{\big\{\sum_{i\ge 1}g(S_i)\mathds{1}_{\{T_i\le t-s\}}\ge(\alpha+\gamma)(t-s)\big\}}U_{t-s}\bigg]\,
  \mathds{1}_{\{\gamma t<s\le\eta t\}}\,p_o(s),
\end{align}
where the fact that $N_{t-s}\le t-s$ has been used in the last
equality.

To continue estimation we need the following lemma, which is proved in appendix \ref{appendix2}.
\begin{lemma}
\label{lemmaux}
$\lim_{t\uparrow\infty}\Ex_o[(N_t/t)U_t]=1/\Ex_o[S_1]^2$.
\end{lemma}
We also observe that
$\lim_{t\uparrow\infty}\Ex_o[\mathds{1}_{\{\sum_{i\ge
      1}g(S_i)\mathds{1}_{\{T_i\le
      \tau\}}\ge(\alpha+\gamma)\tau\}}U_\tau]=0$ since formula
(\ref{mod2}) and proposition \ref{conv} with $g$ in place of $f$ yield
respectively $\Ex_o[\mathds{1}_{\{\sum_{i\ge
      1}g(S_i)\mathds{1}_{\{T_i\le
      \tau\}}\ge(\alpha+\gamma)\tau\}}U_\tau]\le\probab_\tau^c[\sum_{i\ge
    1}g(S_i)\mathds{1}_{\{T_i\le\tau\}}\ge(\alpha+\gamma)\tau]$ and
$\lim_{\tau\uparrow\infty}\probab_\tau^c[\sum_{i\ge
    1}g(S_i)\mathds{1}_{\{T_i\le\tau\}}\ge(\alpha+\gamma)\tau]=0$ as
$\Ex_o[g(S_1)]/\Ex_o[S_1]=1$ and $\alpha+\gamma>1$.  Then, by these
arguments we deduce that there exists $t_2\ge t_1$ such that
$\tau>(1-\eta)t_2$ implies both $\Ex_o[N_\tau U_\tau]\ge
(1-\epsilon)\tau/\Ex_o[S_1]^2$ and $\Ex_o[\mathds{1}_{\{\sum_{i\ge
      1}g(S_i)\mathds{1}_{\{T_i\le
      \tau\}}\ge(\alpha+\gamma)\tau\}}U_\tau]\le
\epsilon/\Ex_o[S_1]^2$.  If $t>t_2$, then the condition $s\le\eta t$
giving $t-s>(1-\eta)t_2$ allows us to replace (\ref{aux2}) with
\begin{align}
  \nonumber
  \mathcal{E}(t)&\ge\frac{1-2\epsilon}{\Ex_o[S_1]^2}\sum_{s\ge 1}(t-s)\,\mathds{1}_{\{\gamma t<s\le\eta t\}}\,p_o(s)-\big[t\,\mathcal{Q}(\gamma t)\big]^2\\
  \nonumber
  &=\frac{1-2\epsilon}{\Ex_o[S_1]^2}\bigg[(1-\gamma)\,t\,\mathcal{Q}(\gamma t)-(1-\eta)\,t\,\mathcal{Q}(\eta t)
  -t\int_\gamma^\eta\mathcal{Q}(xt)dx\bigg]-\big[t\,\mathcal{Q}(\gamma t)\big]^2.
\end{align}
From here, we obtain (\ref{binf1}) by dividing by $t\,\mathcal{Q}(t)$
first and by sending $t$ to infinity later because on the one hand the
limit
$\lim_{t\uparrow\infty}\mathcal{Q}(xt)/\mathcal{Q}(t)=x^{-\kappa}$ is
uniform with respect to $x$ in the compact interval $[\gamma,\eta]$
(see \cite{RV}, theorem 1.5.2), and on the other hand
$\lim_{t\uparrow\infty}t\,\mathcal{Q}(t)=0$ as
$t\,\mathcal{Q}(t)\le\Ex_o[S_1\mathds{1}_{\{S_1>t\}}]$ and
$\Ex_o[S_1]<+\infty$.

\subsection{An upper bound}

Now we show that for all real numbers $\gamma$, $\eta$, and
$\epsilon$ such that $0<\gamma<1-\alpha$, $\gamma<\eta<1$, and
$\epsilon>0$
\begin{equation}
  \limsup_{t\uparrow\infty}\frac{\mathcal{E}(t)}{t\,\mathcal{Q}(t)}\le
  \frac{1+\epsilon}{\Ex_o[S_1]^2}\bigg[(1-\gamma)\gamma^{-\kappa}-(1-\eta)\eta^{-\kappa}-\int_\gamma^\eta \frac{dx}{x^\kappa}\bigg]
  +\eta^{-\kappa}-1.
  \label{primau}
\end{equation}
In the light of (\ref{binf}), this bound leads us to prove
(\ref{mod3}) because sending $\gamma$ to $1-\alpha$, $\eta$ to 1, and
$\epsilon$ to 0 we find
\begin{equation*}
  \limsup_{t\uparrow\infty}\frac{\mathcal{E}(t)}{t\,\mathcal{Q}(t)}\le
  \frac{1}{\Ex_o[S_1]^2}\bigg[\alpha(1-\alpha)^{-\kappa}-\int_{1-\alpha}^1\frac{dx}{x^\kappa}\bigg].
\end{equation*}

Pick $\gamma$, $\eta$, and $\epsilon$ such that $0<\gamma<1-\alpha$,
$\gamma<\eta<1$, and $\epsilon>0$. Obviously, for every $n\ge 1$ and
$t\ge 1$ we have
\begin{equation*}
  1\le\prod_{j=1}^n\mathds{1}_{\{S_j\le\gamma t\}}+\sum_{j=1}^n\mathds{1}_{\{S_j>\gamma t\}}
\end{equation*}
since either $S_j\le\gamma t$ for all $j\le n$ or $S_j>\gamma t$ for
at least one $j\le n$. It follows that for any $t\ge 1$
\begin{equation}
  \mathcal{E}(t)=\sum_{n=1}^t\Ex_o\bigg[\mathds{1}_{\big\{\sum_{i=1}^ng(S_i)\le\alpha t\big\}}\mathds{1}_{\{T_n=t\}}\bigg]\le
  \mathcal{E}_1(t)+\mathcal{E}_2(t)
  \label{auxu1}
\end{equation}
with
\begin{equation*}
\mathcal{E}_1(t):=\sum_{n=1}^t\Ex_o\Bigg[\mathds{1}_{\big\{\sum_{i=1}^ng(S_i)\le\alpha t\big\}}\prod_{j=1}^n\mathds{1}_{\{S_j\le\gamma
      t\}}\mathds{1}_{\{T_n=t\}}\Bigg]
\end{equation*}
and
\begin{equation*}
\mathcal{E}_2(t):=\sum_{n=1}^tn\,\Ex_o\Big[\mathds{1}_{\{S_n>\gamma t\}}\mathds{1}_{\{T_n=t\}}\Big].
\end{equation*}
The term $\mathcal{E}_1(t)$ can be estimated by means of the following
lemma which is proved in appendix \ref{appendix3}. Set $\xi_t:=-\ln
t\,\mathcal{Q}(\gamma t)$ for all $t$ for brevity and notice that
$\lim_{t\uparrow\infty}\xi_t=+\infty$ since
$x\,\mathcal{Q}(x)\le\Ex_o\big[S_1\mathds{1}_{\{S_1>x\}}\big]$ and
$\Ex_o\big[S_1\big]<+\infty$.
\begin{lemma}
  \label{lemmazl}
  For all sufficiently large $t$ there exist two positive real numbers
  $z_t$ and $\lambda_t$ such that $\Ex_o[\mathds{1}_{\{S_1\le\gamma
      t\}}\,e^{z_tS_1-\lambda_tg(S_1)}]\le 1$ and
  $(z_t-\alpha\lambda_t)t\ge\xi_t+\ln|\xi_t|$.
\end{lemma}
Let $z_t$ and $\lambda_t$ be the numbers introduced by lemma
\ref{lemmazl}, so that the inequalities
$0<\Ex_o[\mathds{1}_{\{S_1\le\gamma
    t\}}\,e^{z_tS_1-\lambda_tg(S_1)}]\le 1$ and
$(z_t-\alpha\lambda_t)t\ge\xi_t+\ln|\xi_t|$ hold for all $t>t_1$ with
some $t_1>0$. The fact that $\lambda_t\ge 0$ allows us to invoke the
inequality $\mathds{1}_{\{\sum_{i=1}^ng(S_i)\le\alpha t\}}\le
e^{\alpha\lambda_t t-\lambda_t\sum_{i=1}^ng(S_i)}$ to obtain
\begin{align}
\nonumber
  \mathcal{E}_1(t)&\le e^{\alpha\lambda_tt}\sum_{n=1}^t\Ex_o\Bigg[\prod_{j=1}^n \mathds{1}_{\{S_j\le\gamma
      t\}}\,e^{-\lambda_tg(S_j)}\,\mathds{1}_{\{T_n=t\}}\Bigg]\\
  \nonumber
  &=e^{(\alpha \lambda_t-z_t)t}\sum_{n=1}^t\Ex_o\Bigg[\prod_{j=1}^n \mathds{1}_{\{S_j\le\gamma
      t\}}\,e^{z_tS_j-\lambda_tg(S_j)}\,\mathds{1}_{\{T_n=t\}}\Bigg]\\
  \nonumber
  &=e^{(\alpha \lambda_t-z_t)t}\sum_{n=1}^t\frac{\Ex_o\big[\prod_{j=1}^t \mathds{1}_{\{S_j\le\gamma
      t\}}\,e^{z_tS_j-\lambda_tg(S_j)}\,\mathds{1}_{\{T_n=t\}}\big]}{\Ex_o\big[\mathds{1}_{\{S_1\le\gamma
        t\}}\,e^{z_tS_1-\lambda_tg(S_1)}\big]^{t-n}}.
\end{align}
At this point, we make use of $\Ex_o[\mathds{1}_{\{S_1\le\gamma
    t\}}\,e^{z_tS_1-\lambda_tg(S_1)}]\le 1$ and
$(z_t-\alpha\lambda_t)t\ge\xi_t+\ln|\xi_t|$ to get for any $t>t_1$ at
the upper bound
\begin{align}
\nonumber
  \mathcal{E}_1(t)&\le e^{-\xi_t-\ln|\xi_t|}\sum_{n=1}^t\frac{\Ex_o\big[\prod_{j=1}^t \mathds{1}_{\{S_j\le\gamma
      t\}}\,e^{z_tS_j-\lambda_tg(S_j)}\,\mathds{1}_{\{T_n=t\}}\big]}{\Ex_o\big[\mathds{1}_{\{S_1\le\gamma
        t\}}\,e^{z_tS_1-\lambda_tg(S_1)}\big]^t}\\
  \nonumber
  &=\frac{t\,\mathcal{Q}(\gamma t)}{|\xi_t|}\,\frac{\Ex_o\big[\prod_{j=1}^t \mathds{1}_{\{S_j\le\gamma
      t\}}\,e^{z_tS_j-\lambda_tg(S_j)}\,U_t\big]}{\Ex_o\big[\mathds{1}_{\{S_1\le\gamma
        t\}}\,e^{z_tS_1-\lambda_tg(S_1)}\big]^t}\\
  &\le\frac{t\,\mathcal{Q}(\gamma t)}{|\xi_t|}\,\frac{\Ex_o\big[\prod_{j=1}^t \mathds{1}_{\{S_j\le\gamma
      t\}}\,e^{z_tS_j-\lambda_tg(S_j)}\big]}{\Ex_o\big[\mathds{1}_{\{S_1\le\gamma
        t\}}\,e^{z_tS_1-\lambda_tg(S_1)}\big]^t}=\frac{t\,\mathcal{Q}(\gamma t)}{|\xi_t|}.
\label{uuu1}
\end{align}
As far as $\mathcal{E}_2(t)$ is concerned, we let $\eta$ come into
play by writing for all $t\ge 2$
\begin{align}
\nonumber
\mathcal{E}_2(t)&=\sum_{n=1}^tn\,\Ex_o\Big[\mathds{1}_{\{\gamma t<S_n\le\eta t\}}\mathds{1}_{\{T_n=t\}}\Big]
+\sum_{n=1}^tn\,\Ex_o\Big[\mathds{1}_{\{S_n>\eta t\}}\mathds{1}_{\{T_n=t\}}\Big]\\
\nonumber
&= \sum_{n=2}^t(n-1)\,\Ex_o\Big[\mathds{1}_{\{\gamma t<S_n\le\eta t\}}\mathds{1}_{\{T_{n-1}+S_n=t\}}\Big]
+ \Ex_o\Big[\mathds{1}_{\{\gamma t<S_1\le\eta t\}}U_t\Big]+\Ex_o\Big[\mathds{1}_{\{S_1>\eta t\}}N_tU_t\Big]\\
\nonumber
&=\sum_{s=1}^{t-1}\sum_{n=1}^{t-s}n\,\Ex_o\big[\mathds{1}_{\{T_n=t-s\}}\big]\,\mathds{1}_{\{\gamma t<s\le\eta t\}}\,p_o(s)
+ \Ex_o\big[\mathds{1}_{\{\gamma t<S_1\le\eta t\}}U_t\big]+\Ex_o\big[\mathds{1}_{\{S_1>\eta t\}}N_tU_t\big]\\
\nonumber
&=\sum_{s\ge 1}\Ex_o\big[N_{t-s}U_{t-s}\big]\,\mathds{1}_{\{\gamma t<s\le\eta t\}}\,p_o(s)
+ \Ex_o\big[\mathds{1}_{\{\gamma t<S_1\le\eta t\}}U_t\big]+\Ex_o\big[\mathds{1}_{\{S_1>\eta t\}}N_tU_t\big].
\end{align}
We notice that $\Ex_o[\mathds{1}_{\{\gamma t<S_1\le\eta
    t\}}U_t]\le\Ex_o[\mathds{1}_{\{S_1>\gamma t\}}]=\mathcal{Q}(\gamma
t)$. We also observe that the condition $U_t=1$ implies $S_1\le t$, so
that $\Ex_o[\mathds{1}_{\{S_1>\eta
    t\}}N_tU_t]=\Ex_o[\mathds{1}_{\{\eta t<S_1\le t\}}N_tU_t]$ and
hence $\Ex_o[\mathds{1}_{\{S_1>\eta t\}}N_tU_t]\le
t\,\Ex_o[\mathds{1}_{\{\eta t<S_1\le t\}}]=t\,\mathcal{Q}(\eta
t)-t\,\mathcal{Q}(t)$. Finally, we recall that
$\lim_{\tau\uparrow\infty}\Ex_o[(N_\tau/\tau) U_\tau]=1/\Ex_o[S_1]^2$
by lemma \ref{lemmaux}, which entails that there exists $t_2\ge t_1$
with the property that $\Ex_o[N_\tau U_\tau]\le
(1+\epsilon)\tau/\Ex_o[S_1]^2$ for all $\tau>(1-\eta)t_2$.  This way,
since $t-s\ge(1-\eta)t$ when $s\le\eta t$, for every $t>t_2$ we get
the bound
\begin{align}
\nonumber
\mathcal{E}_2(t)&\le \frac{1+\epsilon}{\Ex_o[S_1]^2}\sum_{s\ge 1}(t-s)\,\mathds{1}_{\{\gamma t<s\le\eta t\}}\,p_o(s)
+\mathcal{Q}(\gamma t)+t\,\mathcal{Q}(\eta t)-t\,\mathcal{Q}(t)\\
\nonumber
&=\frac{1+\epsilon}{\Ex_o[S_1]^2}\bigg[(1-\gamma)\,t\,\mathcal{Q}(\gamma t)-(1-\eta)\,t\,\mathcal{Q}(\eta t)-t\int_\gamma^\eta\mathcal{Q}(xt)dx\bigg]\\
&+\mathcal{Q}(\gamma t)+t\,\mathcal{Q}(\eta t)-t\,\mathcal{Q}(t).
\label{uuu2}
\end{align}

In conclusion, by combining (\ref{auxu1}) with (\ref{uuu1}) and
(\ref{uuu2}) we obtain 
\begin{align}
  \nonumber
  \mathcal{E}(t)&\le\frac{1+\epsilon}{\Ex_o[S_1]^2}\bigg[(1-\gamma)\,t\,\mathcal{Q}(\gamma t)-(1-\eta)\,t\,\mathcal{Q}(\eta t)
    -t\int_\gamma^\eta\mathcal{Q}(xt)dx\bigg]\\
\nonumber
  &+\mathcal{Q}(\gamma t)+\frac{t\,\mathcal{Q}(\gamma t)}{|\xi_t|}+t\,\mathcal{Q}(\eta t)-t\,\mathcal{Q}(t)
\end{align}
for all $t>t_2$. Recalling that the limit
$\lim_{t\uparrow\infty}\mathcal{Q}(xt)/\mathcal{Q}(t)=x^{-\kappa}$ is
uniform with respect to $x$ in the compact interval $[\gamma,\eta]$
and that $\lim_{t\uparrow\infty}\xi_t=+\infty$, (\ref{primau}) follows
from here by dividing by $t\,\mathcal{Q}(t)$ first and by sending $t$
to infinity later.

\section*{Acknowledgments}
The author is grateful to Giambattista Giacomin for a critical reading
of the manuscript and valuable comments.

\appendix

\section{Proof of lemma \ref{lemmaQtail}}
\label{appendix1}

If $\mathcal{L}$ varies slowly and
\begin{equation}
\lim_{x\uparrow+\infty}\frac{\mathcal{Q}(x)}{x^{-\kappa}\mathcal{L}(x)}=\frac{1}{\kappa},
\label{A1}
\end{equation}
then $\mathcal{Q}$ varies regularly with index $-\kappa$.

Let us demonstrate the limit (\ref{A1}).  Pick two real numbers
$\delta\in(0,\kappa)$ and $K>1$ and observe that there exists $x_o>1$
such that $(x/y)^\delta\mathcal{L}(x)/K\le \mathcal{L}(y)\le
K(x/y)^{-\delta}\mathcal{L}(x)$ when $x_o\le x\le y$ (see \cite{RV},
theorem 1.5.6). Recalling that $p_o(s)=s^{-\kappa-1}\mathcal{L}(s)$,
it follows that $x^\delta\mathcal{L}(x) s^{-\kappa-\delta-1}/K\le
p_o(s)\le K x^{-\delta}\mathcal{L}(x)s^{-\kappa+\delta-1}$ if $x_o\le
x<s$. The upper bound gives for $x\ge x_o>1$
\begin{equation*}
  \mathcal{Q}(x)=\sum_{s\ge 1}\mathds{1}_{\{s>x\}}p_o(s)\le \sum_{s\ge 1}\mathds{1}_{\{s>x\}}K x^{-\delta}\mathcal{L}(x)s^{-\kappa+\delta-1}
\le\frac{K}{\kappa-\delta}\frac{x^{-\delta}}{(x-1)^{\kappa-\delta}}\mathcal{L}(x),
\end{equation*}
where the inequality
$(\kappa-\delta)s^{-\kappa+\delta-1}\le(s-1)^{-\kappa+\delta}-s^{-\kappa+\delta}$
valid for all $s>1$ has been used to change the second series with a
telescoping series. Similarly, the lower bound and the inequality
$s^{-\kappa-\delta}-(s+1)^{-\kappa-\delta}\le(\kappa+\delta)s^{-\kappa-\delta-1}$
valid for all $s\ge 1$ yield for $x\ge x_o$
\begin{equation*}
  \mathcal{Q}(x)\ge \sum_{s\ge 1}\mathds{1}_{\{s>x\}}x^\delta\mathcal{L}(x) s^{-\kappa-\delta-1}/K\
\ge\frac{1}{K(\kappa+\delta)}x^{-\kappa}\mathcal{L}(x).
\end{equation*}
This way, we find
\begin{equation*}
  \frac{1}{K(\kappa+\delta)}\le \liminf_{x\uparrow+\infty}\frac{\mathcal{Q}(x)}{x^{-\kappa}\mathcal{L}(x)}
  \le \limsup_{x\uparrow+\infty}\frac{\mathcal{Q}(x)}{x^{-\kappa}\mathcal{L}(x)}\le\frac{K}{\kappa-\delta}
\end{equation*}
and (\ref{A1}) follows from here by sending $\delta$ to 0 and $K$ to
1.

\section{Proof of lemma \ref{lemmaux}}
\label{appendix2}

Recall that $N_t$ is the cumulative reward corresponding to a function
$f$ identically equal to 1.  Since $0\le N_t/t\le 1$ and
$0<1/\Ex_o[S_1]<1$, for any $t\ge 1$ and $\delta>0$ we can write
\begin{align}
  \nonumber
  \Bigg|\Ex_o\bigg[\bigg(\frac{N_t}{t}-\frac{1}{\Ex_o[S_1]}\bigg)U_t\bigg]\Bigg|&\le
  \delta\,\Ex_o\bigg[\mathds{1}_{\big\{\big|\frac{N_t}{t}-\frac{1}{\Ex_o[S_1]}\big|<\delta\big\}}U_t\bigg]
  +2\,\Ex_o\bigg[\mathds{1}_{\big\{\big|\frac{N_t}{t}-\frac{1}{\Ex_o[S_1]}\big|\ge\delta\big\}}U_t\bigg]\\
  \nonumber
  &\le\delta+2\,\probab_t^c\bigg[\bigg|\frac{N_t}{t}-\frac{1}{\Ex_o[S_1]}\bigg|\ge\delta\bigg],
\end{align}
where the second bound is obtained by applying formula (\ref{mod2}) to
$N_t$. From here, we get at the proof of the lemma by combining the
limit $\lim_{t\uparrow\infty}\Ex_o[U_t]=1/\Ex_o[S_1]$ with the
convergence in probability stated by proposition \ref{conv}.

\section{Proof of lemma \ref{lemmazl}}
\label{appendix3}

Recall that $\xi_t:=-\ln t\,\mathcal{Q}(\gamma t)$ and for every $t\ge
1$ set
\begin{equation*}
\lambda_t:=\frac{(1-\gamma)\xi_t-4\kappa\ln|\xi_t|}{\gamma(1-\gamma)\,t}
\end{equation*}
and
\begin{equation*}
z_t:=\alpha\lambda_t+\frac{\xi_t+\ln|\xi_t|}{t}.
\end{equation*}

It is manifest that $(z_t-\alpha\lambda_t)t\ge\xi_t+\ln|\xi_t|$ for
each $t$ and, since $\xi_t$ goes to infinity when $t$ is sent to
infinity, there exists a positive integer $\tau_1$ such that $\xi_t>1$
and the real numbers $\lambda_t$ and $z_t$ are positive for any
$t>\tau_1$. We prove the lemma by showing that for all sufficiently
large $t>\tau_1$
\begin{equation}
\label{lemmazl_1}
  \sum_{s\ge 1}\mathds{1}_{\{s\le\gamma t/\xi_t^2\}}\,e^{z_ts-\lambda_tg(s)}\,p_o(s)\le 1-\frac{2}{t}
\end{equation}
and
\begin{equation}
  \label{lemmazl_2}
\sum_{s\ge 1}\mathds{1}_{\{\gamma t/\xi_t^2<s\le\gamma t\}}\,e^{z_ts-\lambda_tg(s)}\,p_o(s)\le\frac{2}{t}.
\end{equation}
We point out that $\lim_{t\uparrow\infty}\gamma t/\xi_t^2=+\infty$
since $\mathcal{Q}(x)\ge x^{-\kappa-\delta}$ for any $\delta>0$ and
all sufficiently large $x$ as $\mathcal{Q}$ varies regularly with
index $-\kappa$ (see \cite{RV}, proposition 1.3.6).

Let us verify (\ref{lemmazl_1}) at first. Let $M<+\infty$ be a
positive constant such that $|g(s)|\le Ms$ for all $s$, which
certainly exists because $\lim_{s\uparrow\infty}g(s)/s=0$. Then, let
$K<+\infty$ be a positive constant such that
$\gamma(tz_t+Mt\lambda_t)^2\xi_t^{-2}e^{\gamma(tz_t+Mt\lambda_t)\xi_t^{-2}}\le
K$ for all $t>\tau_1$, which exists because
$\lim_{t\uparrow\infty}\gamma(tz_t+Mt\lambda_t)^2\xi_t^{-2}e^{\gamma(tz_t+Mt\lambda_t)\xi_t^{-2}}=\gamma^{-1}(\alpha+\gamma+M)^2$.
Finally, let $\tau_2\ge\tau_1$ be an integer such that $z_ts-\lambda_t
g(s)\ge 0$ for all $t>\tau_2$ and $s>\gamma t/\xi_t^2$, which exists
because $\lim_{t\uparrow\infty}z_t/\lambda_t=\alpha+\gamma>0$,
$\lim_{s\uparrow\infty}g(s)/s=0$, and $\lim_{t\uparrow\infty}\gamma
t/\xi_t^2=+\infty$. Pick $t>\tau_2$. The bound $e^y\le 1+y+y^2e^{|y|}$
valid for all $y\in\Rl$ yields for any $s\le\gamma t/\xi_t^2$
\begin{align}
\nonumber
  e^{z_ts-\lambda_tg(s)}&\le 1+z_ts-\lambda_tg(s)+(z_t+M\lambda_t)^2s^2 e^{(z_t+M\lambda_t)s}\\
  \nonumber
  &\le 1+z_ts-\lambda_tg(s)+\frac{\gamma}{t}\bigg(\frac{tz_t+Mt\lambda_t}{\xi_t}\bigg)^2s\,e^{\gamma\frac{tz_t+Mt\lambda_t}{\xi_t^2}}\\
  \nonumber
  &\le 1+z_ts-\lambda_tg(s)+\frac{K}{t}s.
\end{align}
This bound, in combination with the facts that $z_ts-\lambda_tg(s)\ge
0$ for $s>\gamma t/\xi_t^2$ and that $\Ex_o[g(S_1)]=\Ex_o[S_1]$, gives
\begin{align}
\nonumber
\sum_{s\ge 1}\,\mathds{1}_{\{s\le\gamma t/\xi_t^2\}}\,e^{z_ts-\lambda_tg(s)}\,p_o(s)
&\le 1+\sum_{s\ge 1}\mathds{1}_{\{s\le\gamma t/\xi_t^2\}}\,\big[z_ts-\lambda_tg(s)\big]\,p_o(s)
+\frac{K}{t}\Ex_o[S_1]\\
\nonumber
&\le 1+\Ex_o[S_1]\Bigg\{\frac{tz_t-t\lambda_t}{\xi_t}+\frac{K}{\xi_t}+\frac{2}{\Ex_o[S_1]\xi_t}\Bigg\}\frac{\xi_t}{t}-\frac{2}{t}.
\end{align}
The term between braces goes to $-(1-\alpha-\gamma)/\gamma<0$ when $t$
is sent to infinity, so that it is non-positive for all $t$ larger
than some $\tau_3\ge\tau_2$. It follows that (\ref{lemmazl_1}) holds
for all $t>\tau_3$.

Let us now prove (\ref{lemmazl_2}). To begin with, we observe that for
any $\delta>0$ and $K>1$ there exists $x_o>0$ such that
$y^{\kappa+\delta}\mathcal{Q}(y)\le
K\,x^{\kappa+\delta}\mathcal{Q}(x)$ if $x_o\le y\le x$ (see \cite{RV},
theorem 1.5.6). Thus, taking $\delta=\kappa-\gamma/2$ and $K=2$,
$\delta$ being positive because $\kappa\ge 1$ and $\gamma<1-\alpha$ by
hypothesis, the fact that $\lim_{t\uparrow\infty}\gamma
t/\xi_t^2=+\infty$ implies that an integer $\tau_4\ge\tau_3$ can be
found in such a way that $\mathcal{Q}(\gamma t/\xi_t^2)\le
2\,\xi_t^{4\kappa-\gamma}\mathcal{Q}(\gamma t)$ for all
$t>\tau_4$. Then, since $1-\alpha-\gamma>0$, the limit
$\lim_{s\uparrow\infty}g(s)/s=0$ ensures us that there exists an
integer $\tau_5\ge\tau_4$ with the property that
$g(s)\ge-(1-\alpha-\gamma)s$ for all $s>\gamma t/\xi_t^2$ whenever
$t>\tau_5$. This way, for every $t>\tau_5$ we find the bound
\begin{align}
\nonumber
\sum_{s\ge 1}\mathds{1}_{\{\gamma t/\xi_t^2<s\le\gamma t\}}\,e^{z_ts-\lambda_tg(s)}\,p_o(s)&\le
e^{[z_t+(1-\alpha-\gamma)\lambda_t]\gamma t}\,\mathcal{Q}(\gamma t/\xi_t^2)\\
\nonumber
&=\xi_t^{\gamma-4\kappa}e^{\xi_t}\,\mathcal{Q}(\gamma t/\xi_t^2)=\frac{\xi_t^{\gamma-4\kappa}}{t}\,\frac{\mathcal{Q}(\gamma t/\xi_t^2)}{\mathcal{Q}(\gamma t)}
\le \frac{2}{t}.
\end{align}
In conclusion, we have that both (\ref{lemmazl_1}) and
(\ref{lemmazl_2}) are satisfied if $t>\tau_5$.

\end{document}